\documentclass[11pt]{amsart}
\usepackage{amssymb}
\usepackage[mathscr]{eucal}
\usepackage{hyperref}
\usepackage{url}
\usepackage[all]{xy}
\usepackage{array}
\usepackage{tikz-cd}
\usepackage[T1]{fontenc}
\usepackage[margin=1.15in]{geometry}

\makeatletter
\let\@wraptoccontribs\wraptoccontribs
\makeatother

\newtheorem{thm}[equation]{Theorem}
\newtheorem{lem}[equation]{Lemma}
\newtheorem*{conjecture*}{Conjecture}

\newtheorem{cor}[equation]{Corollary}
\newtheorem{prop}[equation]{Proposition}

\newtheorem{question}{Question}
\theoremstyle{definition}
\newtheorem{rem}[equation]{Remark}
\newtheorem{defn}[equation]{Definition}
\newtheorem{exam}[equation]{Example}

\newcommand{\rank}{\mathrm{rank}}
\numberwithin{equation}{section}
\def\End{\mathrm{End}}

\def\Z{\mathbb{Z}}
\def\Q{\mathbb{Q}}
\def\F{\mathbb{F}}

\def\C{\mathbb{C}}

\def\Fl{\F_\ell}
\def\an{\mathrm{an}}

\def\O{\mathcal{O}}
\def\id{\mathrm{id}}

\def\L{\mathcal{L}}

\def\cS{\mathcal{S}}
\def\T{\mathcal{T}}

\def\Stab{\mathrm{Stab}}

\def\a{\mathfrak{a}}

\def\p{\mathfrak{p}}

\def\Hom{\mathrm{Hom}}
\def\Gal{\mathrm{Gal}}

\def\Aut{\mathrm{Aut}}

\def\tors{\mathrm{tors}}

\def\GL{\mathrm{GL}}

\author{Minseok Kim}
\address{Research Institute of Basic Sciences, Seoul National University, Seoul 08826, Republic of Korea}
\email{\href{mailto:minseokkim.math@gmail.com}{minseokkim.math@gmail.com}}
\author{Myungjun Yu}
\address{Department of Mathematics, Yonsei University, Seoul 03722, Republic of Korea}
\email{\href{mailto:mjyu@yonsei.ac.kr}{mjyu@yonsei.ac.kr}}

\title{Fields Generated by a Single Point on Elliptic Curves}
\begin{document}
\begin{abstract}
We study the finite extensions of a number field generated by a
single point on an elliptic curve. We construct geometrically
non-isomorphic elliptic curves with complex multiplication that
generate the same collection of extensions, giving a negative answer
to a question of Mazur and Rubin. For torsion-generated fields, we
prove that inclusion of the collections up to finitely many
exceptions implies isogeny over the ground field, and in the non-CM
case even isomorphism. When all geometric endomorphisms are defined
over the ground field, we characterize equality of the torsion-generated
collections by the existence, for every prime $p$, of an isogeny over
the ground field of degree prime to $p$. This condition also implies
equality of the full collections of fields generated by a single point.
\end{abstract}
\maketitle

\section{Introduction}

Let $K$ be a number field and fix an algebraic closure $\overline{K}$.
Mazur and Rubin \cite{ds} asked whether a smooth projective curve
is determined, up to geometric isomorphism, by the finite extensions
of $K$ generated by a single point on the curve.
We give a negative answer already for elliptic curves, by constructing
geometrically non-isomorphic curves with complex multiplication (CM)
and the same collection of extensions. For non-CM elliptic curves,
we prove a rigidity result:
the collection generated by torsion points determines the
$K$-isomorphism class, even if one allows finitely many exceptions
in a one-sided inclusion.

We first recall the notion of an extension generated by a point.

\begin{defn}
Let $V$ be a variety over $K$. A finite extension $L/K$ is
\emph{generated by a point of $V$} if $L=K(x)$ for some
$x\in V(\overline K)$. Equivalently, there is a point $x\in V(L)$
that is not defined over any field $L'$ with
$K\subseteq L'\subsetneq L$.
\end{defn}

When $V$ is a variety over $K$, we say that a finite extension
$L/K$ belongs to $V$ if it is generated over $K$ by a point of $V$.
We write
$$
\L(V;K):=\{K(x)/K:x\in V(\overline K)\}
$$
for the collection of all finite extensions of $K$ belonging to $V$,
viewed as subfields of the fixed algebraic closure $\overline K$.

A basic example arises when $V$ contains a curve that is
$K$-isomorphic to an open subset of $\mathbb{P}^1$.
In this case, the primitive element theorem implies that every
finite extension of $K$ belongs to $V$.
Mazur and Rubin prove the following converse for irreducible curves.

\begin{thm}\cite[Theorem 1.8]{ds}
Let $X$ be an irreducible curve over $K$.
Then the following are equivalent:
\begin{enumerate}
    \item all but finitely many finite extensions $L/K$
    belong to $X$;
    \item $X$ is birationally isomorphic over $K$ to the
    projective line.
\end{enumerate}
\end{thm}

This naturally leads to the question of how much information
about a curve is encoded in the collection $\L(X;K)$.
In particular, Mazur and Rubin ask the following question.

\begin{question}\cite[Question 1.9]{ds}\label{question1.9}
Let $X$ and $Y$ be irreducible smooth projective curves over
a number field $K$.
If $\L(X;K)=\L(Y;K)$, are $X$ and $Y$ necessarily
isomorphic over $\overline{K}$?
\end{question}

If $\overline{K}$ is replaced by $K$ in
Question~\ref{question1.9}, the answer is negative;
see \cite[Proposition 2.5]{ds}.

It is also useful to allow the two collections of extensions
to differ by finitely many elements. Thus, we write
$$
\L(X;K)\approx\L(Y;K)
$$
when $\L(X;K)$ and $\L(Y;K)$ differ by at most finitely many
elements. In the genus-zero case, Mazur and Rubin prove
the following. 

\begin{thm}\cite[Theorem 1.10]{ds}
Suppose $X$ and $Y$ are irreducible curves over $K$, and $Y$
has genus zero. Then $\L(X;K)\approx\L(Y;K)$ if and only if
$X$ and $Y$ are birationally isomorphic over $K$.
\end{thm}

To study the rigidity of these collections for elliptic curves,
we first consider extensions generated by a single torsion point.
For an elliptic curve $E/K$, define
$$
\T(E;K)=\{K(P)/K:P\in E(\overline K)_{\mathrm{tors}}\}.
$$
These extensions are controlled by the Galois representations
attached to $E$ and are therefore more accessible than the
full collection $\L(E;K)$.

Related questions have been studied using the full
division fields $K(E[n])$. For example, Lombardo~\cite{lombardo}
studies abelian varieties with identical division fields at every
level. Our collection instead records the fields generated by
a single torsion point, without specifying the order of that point.

For elliptic curves over $K$, we write $E_1\sim_K E_2$ for
$K$-isogeny and $E_1\simeq_K E_2$ for $K$-isomorphism.
Our first result shows that inclusion of the torsion-generated
collections up to finitely many exceptions implies $K$-isogeny.

\begin{thm}\label{mainthm,T}
Let $E_1$ and $E_2$ be elliptic curves over a number field $K$.
If
$$
\T(E_1;K)\setminus\T(E_2;K)
$$
is finite, then $E_1\sim_K E_2$.
\end{thm}

To refine the isogeny conclusion and compare the full collections,
we introduce an invariant of the group of $K$-homomorphisms.

\begin{defn}
Let $E_1$ and $E_2$ be $K$-isogenous elliptic curves.
The \emph{content} of their group of $K$-homomorphisms is
$$
c_K(E_1,E_2):=
\gcd\{\deg\phi:0\neq\phi\in\Hom_K(E_1,E_2)\}.
$$
\end{defn}

\begin{thm}\label{mainthm,content}
Let $E_1$ and $E_2$ be elliptic curves over a number field $K$
that are $K$-isogenous. The following are equivalent:
\begin{enumerate}
    \item $c_K(E_1,E_2)=1$;
    \item for every integer $M\geq1$, there is a $K$-isogeny
    from $E_1$ to $E_2$ of degree relatively prime to $M$.
\end{enumerate}
If these conditions hold, then
$$
\L(E_1;K)=\L(E_2;K).
$$
\end{thm}

For torsion-generated collections, the content criterion also has
a converse under an endomorphism hypothesis. If all geometric
endomorphisms of both curves are defined over $K$,
Theorem~\ref{refine:criterion} shows that
$$
\T(E_1;K)\setminus\T(E_2;K)\text{ is finite}
\quad\Longleftrightarrow\quad
E_1\sim_K E_2\text{ and }c_K(E_1,E_2)=1.
$$
These conditions are also equivalent to
$\T(E_1;K)=\T(E_2;K)$, and each implies
$\L(E_1;K)=\L(E_2;K)$.

For non-CM curves, the endomorphism hypothesis is automatic.
If they are $K$-isogenous, their homomorphism group over $K$
is generated by an isogeny $\phi_0$, and
$$
\deg(n\phi_0)=n^2\deg\phi_0,
\qquad c_K(E_1,E_2)=\deg\phi_0.
$$
Thus $c_K(E_1,E_2)=1$ if and only if $E_1\simeq_K E_2$.
The criterion therefore takes the following form in the non-CM case.

\begin{thm}\label{mainthm,nonCM}
Let $E_1$ and $E_2$ be elliptic curves over a number field $K$.
Suppose neither curve has CM. Then the following are equivalent:
\begin{enumerate}
    \item $\T(E_1;K)\setminus\T(E_2;K)$ is finite;
    \item $\T(E_1;K)=\T(E_2;K)$;
    \item $E_1\simeq_K E_2$.
\end{enumerate}
If these conditions hold, then
$$
\L(E_1;K)=\L(E_2;K).
$$
\end{thm}

The hypothesis concerns only the collections of fields, with no
specified correspondence between torsion points and no requirement
that their orders agree. Even a one-sided inclusion with finitely
many exceptions therefore forces $K$-isomorphism in the non-CM case.

In the CM case, the content condition can hold for geometrically
non-isomorphic curves. Applying Theorem~\ref{mainthm,content} to such
pairs gives the following negative answer to
Question~\ref{question1.9} in genus one.

\begin{thm}\label{mainthm,CM}
For every prime $q\equiv1\pmod4$, there are a number field $K_q$
and elliptic curves $E_{1,q},E_{2,q}$ over $K_q$ such that
$$
\L(E_{1,q};K_q)=\L(E_{2,q};K_q),
\qquad
E_{1,q}\not\simeq_{\overline{K_q}}E_{2,q}.
$$
Moreover, as $q$ varies, all the curves $E_{1,q}$ and $E_{2,q}$
are pairwise non-isomorphic over $\overline{\Q}$.
\end{thm}

The geometric isomorphism classes in Theorem~\ref{mainthm,CM}
vary with $q$. Simultaneous quadratic twists of the pair for $q=5$
also give infinitely many pairs over a fixed number field, distinct
up to $K$-isomorphism (Corollary~\ref{mainthm,quadtwist}).
All these pairs have equal torsion-generated collections as well
(Remark~\ref{torsioncontent}), so the non-CM hypothesis in
Theorem~\ref{mainthm,nonCM} is necessary.

The CM examples leave open the following question for curves
without CM.

\begin{question}\label{question:nonCM}
Let $E_1$ and $E_2$ be elliptic curves over a number field $K$.
Suppose neither curve has CM. If
$$
\L(E_1;K)=\L(E_2;K),
$$
must $E_1$ and $E_2$ be isomorphic over $\overline K$?
\end{question}

We outline the proofs and the organization of the paper.
We first prove the non-CM result directly, and later use it to
complete the criterion of Theorem~\ref{refine:criterion}.
In Section~\ref{nonCMsection}, Serre's open image theorem~\cite{serre}
allows us to recover the residual Galois modules $E_i[\ell]$ from
the torsion-generated collections for every sufficiently large
prime $\ell$. Faltings' theorem~\cite{faltings2} then gives a
$K$-isogeny. To obtain $K$-isomorphism, we must also detect the
integral information at primes dividing the degree of an isogeny.
For a nontrivial cyclic $K$-isogeny, we compare point stabilizers
on the two integral $p$-adic Tate lattices, using Galois elements
that act trivially on all prime-to-$p$ torsion. This produces
infinitely many fields belonging to the first torsion-generated
collection but not the second (Proposition~\ref{separation}).
Applying this to an isogeny of least positive degree proves
Theorem~\ref{mainthm,nonCM}.

Section~\ref{CMsection} completes the isogeny theorem by comparing
split and nonsplit Cartan images in the CM case, and their orbit
degrees with those of $\GL_2(\F_\ell)$ in the mixed case.
Section~\ref{countersection} develops the criterion using the
content $c_K(E_1,E_2)$. For a fixed algebraic point $P$, only
finitely many differences $\sigma(P)-P$ occur. An isogeny with
degree prime to the orders of the nonzero torsion points among
these differences preserves the stabilizer of $P$, hence $K(P)$.
For the converse in Theorem~\ref{refine:criterion}, equality of
point stabilizers forces the CM Tate lattices to be homothetic
at every prime, giving isogenies of degree prime to each prescribed
prime. We then construct the CM pairs using isogenies of degrees
$2$ and $(q+1)/2$.

Section~\ref{ranksection} concludes with a short application:
equality of the full collections forces rank growth over the same
cyclic extensions of any fixed prime degree, up to finitely many
exceptions, after every finite base extension. In particular, the
quadratic twists have positive rank for the same quadratic characters,
up to finitely many exceptions.

We write $G_K=\Gal(\overline K/K)$.
For an elliptic curve $E/K$ and a prime $p$, let $T_pE$
be its $p$-adic Tate module and put
$V_pE=T_pE\otimes_{\Z_p}\Q_p$.
We denote the Galois representation on $E[\ell]$ by
$\rho_{E,\ell}$, using a choice of basis when writing its
image as a subgroup of $\GL_2(\F_\ell)$.
A $K$-isogeny $\phi:E_1\to E_2$ induces a map
$T_p(\phi)\colon T_pE_1\to T_pE_2$ on Tate modules.
After tensoring with $\Q_p$, we obtain a $G_K$-equivariant
isomorphism $V_p(\phi)\colon V_pE_1\to V_pE_2$.

\section{Torsion-generated fields in the non-CM case}
\label{nonCMsection}

We first recover the residual representations, and hence the
$K$-isogeny class, from the torsion-generated collections. A comparison
of integral Tate lattices then strengthens isogeny to isomorphism.
The preliminary lemmas will also be used in the CM case.

\subsection{Preliminary lemmas}
For a variety $X$ over a number field $K$ and a finite extension
$F/K$, write
$$
X_F:=X\times_{K}\mathrm{Spec}(F)
$$
for the base change of $X$ to $F$.

\begin{lem}\label{torsionbasechange}
Let $X,Y$ be varieties over a number field $K$, and let $F/K$ be
a finite extension. Then
$$
\L(X_F;F)\setminus\L(Y_F;F)
\subseteq
\{FM:M\in\L(X;K)\setminus\L(Y;K)\}.
$$
The same inclusion holds with $\L$ replaced by $\T$ when $X$ and
$Y$ are elliptic curves. In particular, finiteness of a one-sided
difference, as well as equality of the collections, persists under
finite base extension.
\end{lem}

\begin{proof}
Let $N=F(P)$ belong to the difference on the left. If
$K(P)\in\L(Y;K)$, choose $Q\in Y(\overline K)$ such that
$K(P)=K(Q)$. Then
$$
N=F(P)=FK(P)=FK(Q)=F(Q),
$$
a contradiction. Thus $K(P)\in\L(X;K)\setminus\L(Y;K)$,
and $N=FK(P)$ has the required form. The same proof applies to
$\T$, taking both $P$ and $Q$ to be torsion points.
\end{proof}

\begin{lem}\label{torsionorbits}
Let $E/F$ be an elliptic curve over a number field, let $\ell\geq3$
be a prime, and put $L=F(E[\ell])$.
\begin{enumerate}
    \item If $\rho_{E,\ell}(G_F)=\GL_2(\Fl)$, then for every
    $0\ne P\in E[\ell]$,
    $$
    [F(P):F]=\ell^2-1.
    $$
    The normal closure of $F(P)/F$ is $L$, and $F(P)/F$ is not
    Galois.
    \item If the image of $\rho_{E, \ell}$ is a full nonsplit Cartan subgroup, then
    $F(P)=L$ and $[L:F]=\ell^2-1$ for every $0\ne P\in E[\ell]$.
    \item  If the image of $\rho_{E, \ell}$ is a full split Cartan subgroup, choose a
    basis in which the image consists of all diagonal matrices.
    For every $P$ with both coordinates nonzero, one has
    $F(P)=L$ and $[L:F]=(\ell-1)^2$.
\end{enumerate}
\end{lem}

\begin{proof}
For (i), the Galois conjugates of $P$ are all the nonzero points
of $E[\ell]$, since $\GL_2(\Fl)$ acts transitively on these points.
Thus the normal closure of $F(P)/F$ is $L$. Choose a basis of $E[\ell]$ whose first vector is $P$. Then $\Gal(L/F(P))$ corresponds to the subgroup
$$
\left\{
\begin{pmatrix}
1&a\\
0&b
\end{pmatrix}
\in \GL_2(\Fl):a\in\Fl,\ b\in\Fl^\times
\right\},
$$
which has order $\ell(\ell-1)$, so (i) follows. 

For (ii), identify $E[\ell]$ with $\F_{\ell^2}$ so that
the Galois image acts by multiplication by elements of
$\F_{\ell^2}^\times$.
This action is free and transitive on the nonzero elements.
Thus the stabilizer of any nonzero point is trivial, giving
$F(P)=L$ and $[L:F]=\ell^2-1$.

For (iii), a diagonal matrix fixing a point $P$ with both
coordinates nonzero is the identity. Therefore, we have $F(P) = L$ and $[L:F]=(\ell-1)^2$.
\end{proof}

\begin{lem}\label{torsionreduction}
Let $E/F$ be an elliptic curve over a number field, and let
$\ell\geq3$ be a prime unramified in $F$. Suppose that $E$ has
good reduction at all primes above $\ell$. Let $N/F$ be a finite
extension whose normal closure $\widetilde N$ contains
$F(\mu_\ell)$. If $N\in\T(E;F)$, then there is a nonzero point
$R\in E[\ell]$ with $F(R)\subseteq N$.

If, moreover, $G_F$ acts transitively on $E[\ell]\setminus\{0\}$,
then $[N:F]\geq\ell^2-1$. In the equality case,
$$
N=F(R),\qquad \widetilde N=F(E[\ell]).
$$
\end{lem}

\begin{proof}
Write $N=F(Q)$ with $Q$ torsion of order $n$. If $\ell\nmid n$,
the N\'eron--Ogg--Shafarevich criterion shows that $F(E[n])/F$
is unramified above $\ell$. This Galois extension contains $N$,
and hence also $\widetilde N$. It would therefore contain
$F(\mu_\ell)$, which is ramified above $\ell$ because $\ell$ is
odd and unramified in $F$. This is impossible.
Thus $\ell\mid n$, and $R=(n/\ell)Q$ has order $\ell$ and is
defined over $N$.

Under the transitivity hypothesis, $[F(R):F]=\ell^2-1$.
This gives the degree bound, with $N=F(R)$ when equality holds.
The conjugates of $R$ are then all the nonzero $\ell$-torsion
points, so the normal closure is $F(E[\ell])$.
\end{proof}

\subsection{Recovering the residual representations}

We now assume that neither $E_1$ nor $E_2$ has CM.

\begin{defn}
A \emph{radial automorphism} of $\GL_2(\F_\ell)$ is an
automorphism of $\GL_2(\F_\ell)$ of the form
$$
A\mapsto(\det A)^mA,
$$
where $m$ is an integer modulo $\ell-1$ satisfying
$\gcd(2m+1,\ell-1)=1$.
\end{defn}

The following two lemmas will be used in the proof of
Proposition~\ref{mod-elliso}.

\begin{lem}\label{directsumofauto}
Let $\ell\geq5$ be a prime.
The group $\Aut(\GL_2(\F_\ell))$ is the direct product of the
group of inner automorphisms and the group of radial automorphisms.
\end{lem}

\begin{proof}
By the classification recalled in \cite[p.~461]{automorphism}, the group $\Aut(\GL_2(\F_\ell))$ is generated by the inner and radial automorphisms, which commute and have trivial intersection.
Alternatively, the assertion follows directly from
\cite[Lemma~3.1]{mayle-wang} with $g=1$, since
$\operatorname{GSp}_2(\Fl)=\GL_2(\Fl)$ and the multiplier is the determinant.
\end{proof}

\begin{lem}\label{innerlem}
Let $\ell\geq5$ be a prime, and for $0\ne v\in\Fl^2$ put
$$
B_v=\{M\in\GL_2(\Fl):Mv=v\}.
$$
If $g\in\Aut(\GL_2(\Fl))$ preserves the determinant and
permutes the set $\{B_v:0\ne v\in\Fl^2\}$, then $g$ is inner.
\end{lem}

\begin{proof}
Every inner automorphism preserves the determinant and permutes the
subgroups $B_v$.
Thus, by Lemma~\ref{directsumofauto}, we may assume that $g$ is radial
and show that $g=\id$. 
Then
$$
g(A)=(\det A)^mA
$$
for some $m$ modulo $\ell-1$ with $\gcd(2m+1,\ell-1)=1$.
Since $g$ preserves the determinant,
$$
\det g(A)=(\det A)^{2m+1}=\det A
$$
for every $A\in\GL_2(\Fl)$. Hence
$$
2m\equiv0\pmod{\ell-1}.
$$
Since $\ell$ is odd, we have
$$
m\equiv0
\quad\text{or}\quad
m\equiv\frac{\ell-1}{2}\pmod{\ell-1}.
$$
The first case gives $g=\id$.
Suppose that $m\equiv(\ell-1)/2$.
Choose a nonsquare $a\in\Fl^\times$ with $a\neq-1$, and put
$$
A=\begin{pmatrix}1&0\\0&a\end{pmatrix}\in B_{e_1}.
$$
Then
$$
g(A)=(\det A)^mA=-A
=\begin{pmatrix}-1&0\\0&-a\end{pmatrix}.
$$
Since $a\neq-1$, the matrix $g(A)$ has no nonzero fixed vector.
Thus $g(A)$ does not belong to $B_v$ for any $0\ne v\in\Fl^2$,
contradicting the assumption that $g$ permutes the set
$\{B_v:0\ne v\in\Fl^2\}$.
Hence $g=\id$.
\end{proof}

\begin{prop}\label{mod-elliso}
Let $E_1$ and $E_2$ be non-CM elliptic curves over a number
field $K$. Suppose that
$$
\T(E_1;K)\setminus\T(E_2;K)
$$
is finite. Then, for every sufficiently large rational prime
$\ell$, we have
$$
K(E_1[\ell])=K(E_2[\ell])
\quad \text{and} \quad E_1[\ell]\simeq E_2[\ell]
\quad\text{as}\quad\Fl[G_K]\text{-modules}.
$$
\end{prop}

\begin{proof}
Put
$$
\cS:=\T(E_1;K)\setminus\T(E_2;K),
$$
and choose $D\geq1$ such that $[F:K]\leq D$ for every
$F\in\cS$. By Serre's open image theorem~\cite{serre}, for every
sufficiently large prime $\ell$ both residual representations are
surjective. Fix such an $\ell\geq5$, large enough that $\ell$ is
unramified in $K$, both curves have good reduction above $\ell$,
and $\ell^2-1>D$.

Let $0\ne P\in E_1[\ell]$ and put $F=K(P)$. By
Lemma~\ref{torsionorbits}(i),
$$
[F:K]=\ell^2-1>D,
$$
and the normal closure of $F/K$ is $K(E_1[\ell])$.
Thus $F\notin\cS$ and $F\in\T(E_2;K)$.
By the Weil pairing, we have $K(\mu_\ell) \subseteq K(E_1[\ell])$.
Since the residual image of $E_2$ is transitive on nonzero vectors,
the equality case of Lemma~\ref{torsionreduction} gives a point
$0\ne R\in E_2[\ell]$ such that
$$
K(R)=F,
\qquad K(E_2[\ell])=K(E_1[\ell]).
$$

Put $L=K(E_1[\ell])=K(E_2[\ell])$. After choosing bases, the two
residual representations induce isomorphisms
$$
\rho_{E_i,\ell}:\Gal(L/K)\xrightarrow{\sim}\GL_2(\Fl).
$$
They therefore define an automorphism
$$
g=\rho_{E_2,\ell}\circ\rho_{E_1,\ell}^{-1}
\in\Aut(\GL_2(\Fl)).
$$
The determinants of both representations are the mod-$\ell$
cyclotomic character, so $g$ preserves the determinant.
The preceding argument applies to every nonzero $P\in E_1[\ell]$
and gives a nonzero $R\in E_2[\ell]$ with $K(P)=K(R)$.
Writing $B_P,B_R$ for the corresponding vector stabilizers, we get
$$
\begin{aligned}
g(B_P)
&=\rho_{E_2,\ell}\bigl(\Gal(L/K(P))\bigr)\\
&=\rho_{E_2,\ell}\bigl(\Gal(L/K(R))\bigr)=B_R.
\end{aligned}
$$
Thus $g$ maps the finite set of nonzero-vector stabilizers into
itself. This map is injective, and hence a permutation.
Lemma~\ref{innerlem} implies that $g$ is inner, so
$$
E_1[\ell]\simeq E_2[\ell]
\quad\text{as}\quad\Fl[G_K]\text{-modules}.
$$
\end{proof}

\begin{prop}\label{isogprop,nonCM}
Let $E_1$ and $E_2$ be non-CM elliptic curves over a number
field $K$. If
$$
\T(E_1;K)\setminus\T(E_2;K)
$$
is finite, then $E_1$ and $E_2$ are $K$-isogenous.
\end{prop}

\begin{proof}
The result follows from Proposition~\ref{mod-elliso} and the implication
(iii)$\Rightarrow$(i) in \cite[Corollary~2]{faltings2}.
\end{proof}

\subsection{From isogeny to isomorphism}

The following lemma records the forms of adelic openness used here
and in Section~\ref{countersection}. Its non-CM case will strengthen
$K$-isogeny to $K$-isomorphism.

\begin{lem}\label{adelic}
Let $E/K$ be an elliptic curve with
$\End_K(E)=\End_{\overline K}(E)$, and fix a rational prime $p$.
For every sufficiently large integer $N$, the following hold.
\begin{enumerate}
    \item If $E$ is non-CM, fix a $\Z_p$-basis of $T_pE$.
    For every $A\in1+p^NM_2(\Z_p)$, there is $\sigma\in G_K$
    acting on $T_pE$ by $A$ and on every $T_qE$, $q\ne p$,
    by the identity.
    \item If $E$ has CM by an order in $M$, identify $V_pE$
    with the rank-one module $M\otimes\Q_p$ and put
    $R=\O_M\otimes\Z_p$. For every $a\in R$, there is
    $\sigma\in G_K$ acting on $V_pE$ by multiplication by
    $1+p^Na$ and on every $T_qE$, $q\ne p$, by the identity.
\end{enumerate}
In either case, $\sigma$ also acts trivially on $T_qE'$ for
every $q\ne p$ and every elliptic curve $E'$ that is
$K$-isogenous to $E$.
\end{lem}

\begin{proof}
Choose a $\mathbb Z_q$-basis of $T_qE$ for each prime $q$, and let
$$
\rho_{\mathrm{ad}}\colon G_K \longrightarrow
\prod_q \mathrm{GL}_2(\mathbb Z_q)
= \mathrm{GL}_2(\widehat{\mathbb Z})
$$
denote the resulting adelic representation.
If $E$ is non-CM, Serre's open image theorem~\cite{serre} shows that
$\rho_{\mathrm{ad}}(G_K)$ is open in $\GL_2(\widehat{\Z})$.
If $E$ has CM, let $\mathcal O$ be the CM order. By
\cite[\S4.5, Corollaire]{serre}, the image is open in
$(\mathcal O\otimes\widehat{\Z})^\times$, which acts by multiplication
on the Tate modules.
In either case, the image contains the kernel of reduction modulo
some positive integer $m$ in the corresponding group.

For (i), write $m=p^r m_0$ with $p\nmid m_0$, and let $N\geq\max\{1,r\}$.
Given $A\in I+p^N M_2(\mathbb Z_p)$, define
$g=(g_q)_q\in\prod_q\mathrm{GL}_2(\mathbb Z_q)$ by
$$
g_p=A,\qquad g_q=I\quad(q\ne p).
$$
Since $m_0$ is a unit in $\mathbb Z_p$, we have
$m\mathbb Z_p=p^r\mathbb Z_p$. Thus $N\geq r$ gives
$$
g_p-I=A-I\in p^N M_2(\mathbb Z_p)
\subseteq mM_2(\mathbb Z_p).
$$
For every $q\ne p$, we have $g_q-I=0$.
Hence $g$ lies in the kernel of reduction modulo $m$, which is
contained in $\rho_{\mathrm{ad}}(G_K)$.
Choose $\sigma\in G_K$ with $\rho_{\mathrm{ad}}(\sigma)=g$.
Then $\sigma$ acts on $T_pE$ by $A$ and on every $T_qE$ with
$q\ne p$ by the identity.

For (ii), put $\mathcal O_p=\mathcal O\otimes\mathbb Z_p$.
Since $\mathcal O_p$ has finite index in $R$, we have
$p^N R\subseteq m\mathcal O_p\cap p\mathcal O_p$ for all sufficiently
large $N$. Thus, for every $a\in R$, the element $1+p^Na$ is a unit
of $\mathcal O_p$ congruent to $1$ modulo $m\mathcal O_p$.
The tuple with $p$-component $1+p^Na$ and all other components equal
to $1$ therefore lies in the kernel of reduction modulo $m$, and hence
in $\rho_{\mathrm{ad}}(G_K)$.

Finally, a $K$-isogeny gives a $G_K$-equivariant isomorphism
$V_qE\simeq V_qE'$. An element acting trivially on $T_qE$ acts
trivially on both rational Tate modules, and hence on $T_qE'$.
\end{proof}

\begin{prop}\label{separation}
Let $E_1,E_2/K$ be non-CM elliptic curves and let
$\phi:E_1\to E_2$ be a cyclic $K$-isogeny of degree $d>1$.
For every prime $p$ dividing $d$, there are infinitely many
distinct fields $K(P)$, with $P\in E_1[p^\infty]\setminus\{0\}$,
such that
$$
K(P)\notin\T(E_2;K).
$$
\end{prop}

\begin{proof}
Fix a prime $p$ dividing $d$, and put $r=v_p(d)$.
Via $V_p(\phi)^{-1}$, view $T_pE_2$ as a lattice in $V_pE_1$
containing $T_pE_1$.
Since the $p$-primary part of $\ker(\phi)$ is cyclic of order
$p^r$, we may choose a basis $e_1,e_2$ of $T_pE_1$ such that
$$
T_pE_1=\Z_pe_1\oplus\Z_pe_2,
\qquad
T_pE_2=\Z_p(e_1/p^r)\oplus\Z_pe_2.
$$
Write these lattices as $\Lambda_1$ and $\Lambda_2$, respectively,
and put $V=V_pE_1$.
Via the above identification, we identify $E_i[p^\infty]$
with $V/\Lambda_i$ for $i=1,2$.

Let $N\geq r+1$ be any sufficiently large integer such that
Lemma~\ref{adelic} applies to $E_1$. Put
$$
P_N=e_1/p^{N+1}\bmod\Lambda_1.
$$
We claim that
$$
K(P_N)\ne K(Q)
\qquad\text{for every }Q\in E_2(\overline K)_{\tors}.
$$
Suppose that $K(P_N)=K(Q)$ for some torsion point
$Q\in E_2(\overline K)$.
Write $Q=Q_p+Q_{p'}$, where $Q_p$ has $p$-power order
and $Q_{p'}$ has order prime to $p$.
We may write
$$
Q_p=\alpha e_1+\beta e_2\bmod\Lambda_2,
\qquad
\alpha,\beta\in\Q_p.
$$

For $x\in\Z_p$, put
$$
D(x)=\begin{pmatrix}1+p^Nx&0\\0&1\end{pmatrix},
\qquad
U(x)=\begin{pmatrix}1&0\\p^Nx&1\end{pmatrix}.
$$
By Lemma~\ref{adelic}(i), there are
$\sigma_D(x),\sigma_U(x)\in G_K$ whose actions on $T_pE_1$
are $D(x),U(x)$, respectively, and whose actions on $T_qE_1$
are the identity for every $q\ne p$.
The final assertion of Lemma~\ref{adelic} shows that these elements
also act trivially on $T_qE_2$ for $q\ne p$, so they fix $Q_{p'}$.

We have
$$
(D(x)-1)P_N=(x/p)e_1\bmod\Lambda_1,
\qquad
(U(x)-1)P_N=(x/p)e_2\bmod\Lambda_1.
$$
Thus
$$
\sigma_D(x)(P_N)=P_N
\quad\Longleftrightarrow\quad x\in p\Z_p,
$$
and
$$
\sigma_U(x)(P_N)=P_N
\quad\Longleftrightarrow\quad x\in p\Z_p.
$$
On the other hand,
$$
\sigma_D(x)(Q)=Q
\quad\Longleftrightarrow\quad
p^Nx\alpha\in p^{-r}\Z_p,
$$
and
$$
\sigma_U(x)(Q)=Q
\quad\Longleftrightarrow\quad
p^Nx\alpha\in\Z_p.
$$
Since $K(P_N)=K(Q)$, their stabilizers in $G_K$ agree.
Hence, for every $x\in\Z_p$,
\begin{align*}
x\in p\Z_p &\quad\Longleftrightarrow\quad p^Nx\alpha\in p^{-r}\Z_p, \\
x\in p\Z_p &\quad\Longleftrightarrow\quad p^Nx\alpha\in\Z_p.
\end{align*}
Taking $x=p$ in the second equivalence gives $p^{N+1}\alpha\in\Z_p$.
Since $r\geq 1$, it follows that
$p^N\alpha\in p^{-1}\Z_p\subseteq p^{-r}\Z_p$.
The first equivalence with $x=1$ then gives $1\in p\Z_p$, a contradiction.
Therefore,
$$
K(P_N)\notin\T(E_2;K).
$$

This construction works for every sufficiently large $N$,
and $P_N$ has order $p^{N+1}$.
If only finitely many fields $K(P_N)$ occurred, their compositum
would be a number field $F$ over which all these points
are defined.
This would contradict the finiteness of $E_1(F)_{\tors}$.
Hence infinitely many distinct fields $K(P_N)$ occur.
\end{proof}

We now prove Theorem~\ref{mainthm,nonCM}.

\begin{proof}[Proof of Theorem~\ref{mainthm,nonCM}]
If $E_1\simeq_K E_2$, then the fields generated by torsion
points agree.
Thus (iii) implies (ii), and (ii) clearly implies (i).

Suppose that $\T(E_1;K)\setminus\T(E_2;K)$ is finite.
By Proposition~\ref{isogprop,nonCM}, the curves are $K$-isogenous.
Choose a $K$-isogeny $\phi:E_1\to E_2$ of least positive degree.

Suppose that $\ker(\phi)$ is not cyclic.
Then there is a prime $p$ such that $E_1[p]\subseteq\ker(\phi)$.
By \cite[Corollary~III.4.11]{silverman1}, there is a unique isogeny
$\psi:E_1\to E_2$ such that $\phi=\psi\circ[p]$.
Since $\phi$ and $[p]$ are defined over $K$, the uniqueness
of $\psi$ shows that $\psi$ is also defined over $K$. But
$$
\deg\psi=\deg\phi/p^2<\deg\phi,
$$
contradicting the choice of $\phi$.
Hence $\ker(\phi)$ is cyclic.

If $\deg\phi>1$, Proposition~\ref{separation} shows that
$\T(E_1;K)\setminus\T(E_2;K)$ is infinite, a contradiction.
Therefore $\deg\phi=1$, and $E_1\simeq_K E_2$.
This proves that (i) implies (iii).
\end{proof}

\section{The CM and mixed cases}\label{CMsection}

In this section, we assume that at least one of the two
elliptic curves has CM.
We prove that if both curves have CM and
$\T(E_1;K)\setminus\T(E_2;K)$ is finite, then $E_1\sim_K E_2$.
If exactly one curve has CM, both differences of the
torsion-generated collections are infinite.

Let $E/F$ be an elliptic curve with CM by an order $\O$, and suppose
that all endomorphisms of $E$ are defined over $F$.
By the CM open image theorem used in Lemma~\ref{adelic}, for every
sufficiently large prime $\ell$, the image $\rho_{E,\ell}(G_F)$
of the mod-$\ell$ Galois representation is the full group
$(\O/\ell\O)^\times$ acting by multiplication on
$E[\ell]\simeq\O/\ell\O$.
This is a split or nonsplit Cartan subgroup according as $\ell$
splits or is inert in the CM field.
The same conclusion holds with $G_F$ replaced by $G_{F'}$ for
any fixed finite extension $F'/F$ and all sufficiently large $\ell$.

\begin{prop}\label{cmisogthm}
Let $E_1$ and $E_2$ be elliptic curves over a number field $K$,
and assume that both curves have CM.
If
$$
\T(E_1;K)\setminus\T(E_2;K)
$$
is finite, then $E_1$ and $E_2$ are $K$-isogenous.
\end{prop}

\begin{proof}
Let
$$
M_i:=\End_{\overline K}(E_i)\otimes\Q
\qquad(i=1,2),
$$
where we identify each $M_i$ with an imaginary quadratic
subfield of $\overline K$, and put $F=KM_1M_2$.
All the endomorphisms of both curves are defined over $F$.
By Lemma~\ref{torsionbasechange}, the set
$$
\cS_F:=\T(E_1;F)\setminus\T(E_2;F)
$$
is finite.
Choose an integer $D\geq1$ such that $[N:F]\leq D$
for every $N\in\cS_F$.

We first claim that $M_1=M_2$.
Suppose that $M_1\ne M_2$.
By the Chebotarev density theorem, there are infinitely many
rational primes $\ell$ which split in $M_1$ and are inert in $M_2$.
Choose such an $\ell$ sufficiently large that $\ell$ is
unramified in $F$, both curves have good reduction at all
primes above $\ell$, the residual images are the full
Cartan subgroups, and
$$
(\ell-1)^2>D.
$$
The condition on the residual images holds for all sufficiently
large $\ell$ by the CM open image theorem.

Choose $P\in E_1[\ell]$ with coordinates $(1,1)$ in a split
Cartan basis, and put $N=F(P)$. Lemma~\ref{torsionorbits}(iii)
gives
$$
N=F(E_1[\ell]),\qquad [N:F]=(\ell-1)^2>D.
$$
Thus $N\notin\cS_F$ and $N\in\T(E_2;F)$.
The Weil pairing gives $F(\mu_\ell)\subseteq N$.
The full nonsplit Cartan image of $E_2$ is transitive on its
nonzero $\ell$-torsion points, so Lemma~\ref{torsionreduction}
implies
$$
\ell^2-1\leq[N:F]=(\ell-1)^2,
$$
a contradiction.

Thus $M_1=M_2=:M$, and hence $F=KM$.
Since two CM elliptic curves with the same CM field are
isogenous over $\overline K$, we may choose an isogeny
$$
\phi:E_1\longrightarrow E_2
$$
defined over a finite extension $L/F$.

We now claim that $\phi$ is defined over $F$.
Let $\ell$ be a sufficiently large rational prime which
is inert in $M$.
Applying the CM open image theorem over both $F$ and the fixed
extension $L$, the images of the residual representations
of both curves over these fields are the full nonsplit
Cartan subgroups.
In particular,
$$
\rho_{E_i,\ell}(G_L)=\rho_{E_i,\ell}(G_F)
\cong\F_{\ell^2}^\times
\qquad(i=1,2).
$$
We also assume that $\ell$ is unramified in $F$, both curves
have good reduction at all primes above $\ell$,
$\ell\nmid\deg\phi$, and $\ell^2-1>D$.

Choose $0\ne P\in E_1[\ell]$ and put $N=F(P)$.
By Lemma~\ref{torsionorbits}(ii),
$$
N=F(E_1[\ell]),\qquad [N:F]=\ell^2-1>D.
$$
Hence $N\notin\cS_F$ and $N\in\T(E_2;F)$.
Since $N$ is Galois over $F$ and contains $F(\mu_\ell)$,
the equality case of Lemma~\ref{torsionreduction} gives
$$
N=F(E_1[\ell])=F(E_2[\ell]).
$$

Since $\ell\nmid\deg\phi$, the restriction of $\phi$ gives
a $G_L$-equivariant isomorphism
$$
\phi_\ell:E_1[\ell]\xrightarrow{\sim}E_2[\ell].
$$
Both $G_F$-actions factor through $\Gal(N/F)$.
Moreover,
$$
\rho_{E_1,\ell}(G_L)=\rho_{E_1,\ell}(G_F)
$$
shows that the restriction map
$$
G_L\longrightarrow\Gal(N/F)
$$
is surjective.
Thus the $G_L$-equivariance of $\phi_\ell$ implies
its $G_F$-equivariance.

This holds for every sufficiently large prime $\ell$
which is inert in $M$.
Consequently, for each $\sigma\in G_F$,
$$
(\sigma\phi-\phi)|_{E_1[\ell]}=0
$$
for infinitely many $\ell$.
A nonzero homomorphism between elliptic curves has finite
kernel, so $\sigma\phi-\phi=0$.
Therefore $\phi$ is defined over $F$.

If $F=K$, the proof is complete.
Suppose that $F\ne K$, and let $\tau$ be the nontrivial
element of $\Gal(F/K)$.
If
$$
\phi+\tau\phi\ne0,
$$
then this sum is a nonzero homomorphism defined over $K$,
and hence a $K$-isogeny.

Suppose instead that $\tau\phi=-\phi$.
Choose a nonzero endomorphism
$\alpha\in\End_F(E_2)$ whose trace in $M/\Q$ is zero.
Since $\tau$ acts on $\End_F(E_2)$ by complex conjugation,
we have $\tau\alpha=-\alpha$. Thus
$$
\tau(\alpha\circ\phi)
=(\tau\alpha)\circ(\tau\phi)
=(-\alpha)\circ(-\phi)
=\alpha\circ\phi.
$$
Hence $\alpha\circ\phi$ is a $K$-isogeny.
This proves that $E_1$ and $E_2$ are $K$-isogenous.
\end{proof}

\begin{prop}\label{mixedcase}
Let $E_1$ and $E_2$ be elliptic curves over a number field $K$.
Assume that $E_1$ has CM and $E_2$ does not have CM.
Then both
$$
\T(E_1;K)\setminus\T(E_2;K)
\quad\text{and}\quad
\T(E_2;K)\setminus\T(E_1;K)
$$
are infinite.
\end{prop}

\begin{proof}
Let $M=\End_{\overline K}(E_1)\otimes\Q$, identify $M$ with an
imaginary quadratic subfield of $\overline K$, and put $F=KM$.
By Lemma~\ref{torsionbasechange}, it suffices to prove both
assertions over $F$.

All endomorphisms of $E_1$ are defined over $F$, so the CM theory
shows that $F(E_1[n])/F$ is abelian for every $n\geq1$.
Consequently, every field in $\T(E_1;F)$ is Galois over $F$.
For every sufficiently large prime $\ell$, Serre's open image
theorem gives $\rho_{E_2,\ell}(G_F)=\GL_2(\Fl)$.
Choose $0\ne P_\ell\in E_2[\ell]$.
By Lemma~\ref{torsionorbits}(i), $F(P_\ell)/F$ is not Galois
and has degree $\ell^2-1$. Hence
$$
F(P_\ell)\in\T(E_2;F)\setminus\T(E_1;F),
$$
and the unbounded degrees give infinitely many such fields.

For the other direction, let $\ell$ vary through sufficiently
large primes split in $M$. We may assume that $\ell$ is
unramified in $F$, both curves have good reduction above $\ell$,
the image for $E_1$ is the full split Cartan, and the image for
$E_2$ is $\GL_2(\Fl)$. Choose $P_\ell\in E_1[\ell]$ with
coordinates $(1,1)$ in a split Cartan basis. By
Lemma~\ref{torsionorbits}(iii), the field $N_\ell=F(P_\ell)$
satisfies
$$
N_\ell=F(E_1[\ell]),\qquad [N_\ell:F]=(\ell-1)^2.
$$
In particular, $F(\mu_\ell)\subseteq N_\ell$.
If $N_\ell\in\T(E_2;F)$, Lemma~\ref{torsionreduction} would
give $\ell^2-1\leq[N_\ell:F]=(\ell-1)^2$, a contradiction.
Thus
$$
N_\ell\in\T(E_1;F)\setminus\T(E_2;F).
$$
Again the degrees are unbounded, so this difference is infinite.
\end{proof}

We now prove Theorem~\ref{mainthm,T}.

\begin{proof}[Proof of Theorem~\ref{mainthm,T}]
Suppose that $\T(E_1;K)\setminus\T(E_2;K)$ is finite.
By Proposition~\ref{mixedcase}, either both curves have CM
or neither does.
Then the result follows from
Propositions~\ref{cmisogthm} and~\ref{isogprop,nonCM}.
\end{proof}

\section{Equality criteria and CM counterexamples}
\label{countersection}

We now pass from torsion-generated fields to fields generated by
arbitrary points.
We first prove Theorem~\ref{mainthm,content}, which shows that
$\L(E_1;K)=\L(E_2;K)$ whenever $E_1$ and $E_2$ are $K$-isogenous
and $c_K(E_1,E_2)=1$. We then characterize when
$\T(E_1;K)=\T(E_2;K)$ under the assumption that all endomorphisms
are defined over $K$. Finally, we apply Theorem~\ref{mainthm,content}
to construct counterexamples to Question~\ref{question1.9}.

\subsection{Isogenies of relatively prime degrees}

For a given point $P\in E_1(\overline K)$, the following lemma
gives a condition on the degree of a $K$-isogeny $\phi:E_1\to E_2$
that ensures $K(\phi(P))=K(P)$.
For $P\in E(\overline K)$, let
$$
\Stab_K(P):=\{\sigma\in G_K:\sigma(P)=P\},
$$
so that $K(P)$ is the fixed field of $\Stab_K(P)$.

\begin{lem}\label{isoglem}
Let $E_1$ and $E_2$ be elliptic curves over a number field $K$
and let $P\in E_1(\overline K)$.
There exists an integer $m=m(P)\geq1$ with the following property:
if
$$
\phi:E_1\to E_2
$$
is a $K$-isogeny satisfying $\gcd(\deg\phi,m)=1$, then
$$
K(\phi(P))=K(P).
$$
\end{lem}

\begin{proof}
The set $\Delta(P)=\{\sigma(P)-P:\sigma\in G_K\}$ is finite.
Let $m(P)$ be the least common multiple of the orders of the
nonzero torsion points in $\Delta(P)$, or $1$ if there are none.

If $\gcd(\deg\phi,m(P))=1$, every point of
$\Delta(P)\cap\ker\phi$ has order dividing both integers, so
$\Delta(P)\cap\ker\phi=\{0\}$. Since $\phi$ is defined over $K$,
$$
\sigma(\phi(P))=\phi(P)
\quad\Longleftrightarrow\quad
\sigma(P)=P.
$$
Thus $\Stab_K(\phi(P))=\Stab_K(P)$, and
$K(\phi(P))=K(P)$.
\end{proof}

\begin{prop}\label{isoprop}
Let $E_1$ and $E_2$ be elliptic curves over a number field $K$.
Suppose that for every integer $M\geq1$, there exists
a $K$-isogeny
$$
\phi_M:E_1\to E_2
\qquad\text{with}\qquad
\gcd(\deg\phi_M,M)=1.
$$
Then $\L(E_1;K)=\L(E_2;K)$.
\end{prop}

\begin{proof}
For $P\in E_1(\overline K)$, apply the hypothesis with $M=m(P)$
from Lemma~\ref{isoglem}. Then $K(\phi_M(P))=K(P)$, proving
$\L(E_1;K)\subseteq\L(E_2;K)$. The duals of the isogenies
$\phi_M$ have the same degrees and are defined over $K$, so
the reverse inclusion follows by the same argument.
\end{proof}

We now prove Theorem~\ref{mainthm,content}, which we restate here.

\begin{thm}[=Theorem~\ref{mainthm,content}]
\label{degreecontentprop}
Let $E_1$ and $E_2$ be elliptic curves over a number field $K$
that are $K$-isogenous.
Then the following conditions are equivalent:
\begin{enumerate}
    \item $c_K(E_1,E_2)=1$;\label{H1}
    \item for every integer $M\geq1$, there is a $K$-isogeny
    from $E_1$ to $E_2$ of degree relatively prime to $M$.\label{H2}
\end{enumerate}
If these conditions hold, then
$$
\L(E_1;K)=\L(E_2;K).
$$
\end{thm}

\begin{proof}
Condition~\ref{H2}, applied with $M=p$ for each prime $p$,
implies that no prime divides every isogeny degree. Thus it
implies condition~\ref{H1}.

Conversely, suppose $c_K(E_1,E_2)=1$ and fix $M>1$.
For each prime $p\mid M$, choose a $K$-isogeny $\psi_p$ whose
degree is prime to $p$. By \cite[Corollary III.7.5]{silverman1}, the group $\Hom_K(E_1,E_2)$ is a free
$\Z$-module of finite rank, so the Chinese remainder theorem gives
$\phi\in\Hom_K(E_1,E_2)$ satisfying
$$
\phi\equiv\psi_p\pmod{p\Hom_K(E_1,E_2)}
\qquad(p\mid M).
$$
With the convention $\deg0=0$, the degree map is an integral
quadratic form on $\Hom_K(E_1,E_2)$ by \cite[Corollary III.6.3]{silverman1}. Hence
$\deg\phi\equiv\deg\psi_p\not\equiv0\pmod p$ for every
$p\mid M$. Thus $\phi\ne0$ and $\gcd(\deg\phi,M)=1$.
The case $M=1$ follows from the existence of any nonzero
$K$-homomorphism. This proves condition~\ref{H2}, and
Proposition~\ref{isoprop} gives equality of the collections $\L$.
\end{proof}

\begin{cor}\label{coprimedegreecor}
Let $E_1$ and $E_2$ be elliptic curves over a number field $K$.
Suppose that there exist $K$-isogenies
$$
\phi,\psi:E_1\longrightarrow E_2
\qquad\text{such that}\qquad
\gcd(\deg\phi,\deg\psi)=1.
$$
Then $\L(E_1;K)=\L(E_2;K)$.
\end{cor}

\begin{proof}
Since $c_K(E_1,E_2)$ divides both $\deg\phi$ and $\deg\psi$,
we have $c_K(E_1,E_2)=1$.
The result follows from Theorem~\ref{degreecontentprop}.
\end{proof}

\begin{rem}\label{torsioncontent}
If $c_K(E_1,E_2)=1$, then $\T(E_1;K)=\T(E_2;K)$ as well.
Indeed, the proof of Proposition~\ref{isoprop} applies to
torsion points, since isogenies send torsion points to
torsion points.
In the proof of Theorem~\ref{mainthm,CM}, we construct geometrically
non-isomorphic CM elliptic curves $E_1$ and $E_2$ over a number field
$K$ with $c_K(E_1,E_2)=1$. Thus $\T(E_1;K)=\T(E_2;K)$ for these
curves, showing that the non-CM hypothesis in
Theorem~\ref{mainthm,nonCM} is necessary.
\end{rem}

\subsection{Endomorphisms defined over the ground field}

We next prove that, when all endomorphisms are defined over $K$,
the content condition is also necessary for equality of the
torsion-generated collections. We use the following lattice observation.

\begin{lem}\label{refine:smallsubmodule}
Let $V$ be a finite-dimensional $\Q_p$-vector space, let $R$ be a
full $\Z_p$-lattice in $V$, and let $J\subseteq V$ be a
$\Z_p$-submodule. If $J\cap R\subseteq pR$, then $J\subseteq pR$.
\end{lem}

\begin{proof}
If $x\in J\setminus pR$, there is an integer $k\geq0$ such that
$p^kx\in R\setminus pR$. Indeed, take $k=0$ if $x\in R$, and
otherwise take the least $k$ for which $p^kx\in R$. Since $J$ is
a $\Z_p$-submodule, this contradicts $J\cap R\subseteq pR$.
\end{proof}

\begin{thm}\label{refine:criterion}
Let $E_1,E_2/K$ be elliptic curves over a number field $K$, and suppose
that
$$
\End_K(E_i)=\End_{\overline K}(E_i)\qquad(i=1,2).
$$
Then the following conditions are equivalent:
\begin{enumerate}
    \item $\T(E_1;K)\setminus\T(E_2;K)$ is finite;
    \item $\T(E_1;K)=\T(E_2;K)$;
    \item $E_1\sim_K E_2$ and $c_K(E_1,E_2)=1$.
\end{enumerate}
If these conditions hold, then
$$
\L(E_1;K)=\L(E_2;K).
$$
\end{thm}

\begin{proof}
The implication (ii)$\Rightarrow$(i) is immediate. If (iii) holds,
Theorem~\ref{mainthm,content} gives equality of the collections $\L$.
Remark~\ref{torsioncontent} also gives equality of the collections
$\T$.
Thus (iii) implies (ii) and the final assertion.

Suppose that (i) holds. By Theorem~\ref{mainthm,T}, the curves are
$K$-isogenous, so either both curves have CM or neither does.
If neither curve has CM, Theorem~\ref{mainthm,nonCM} gives
$E_1\simeq_K E_2$, and hence $c_K(E_1,E_2)=1$.
It remains to treat the case when both curves have CM.

Choose a $K$-isogeny $\phi:E_1\to E_2$. It induces an isomorphism
$$
\End_K(E_1)\otimes_{\Z}\Q
\xrightarrow{\ \sim\ }
\End_K(E_2)\otimes_{\Z}\Q,
\qquad
\alpha\longmapsto
\frac{1}{\deg\phi}\,\phi\circ\alpha\circ\widehat{\phi},
$$
where $\widehat{\phi}$ is the dual isogeny. Via this isomorphism,
we identify the two algebras and denote the common imaginary
quadratic field by $M$.

Fix a rational prime $p$, and put
$$
A=M\otimes_{\Q}\Q_p,
\qquad R=\O_M\otimes_{\Z}\Z_p.
$$
By \cite[\S4.5, p.~301]{serre}, $V_pE_1$ is free of rank one over $A$. 
Let $V=V_pE_1$ and identify $V_pE_2$ with $V$ via $V_p(\phi)^{-1}$.
Choose an $A$-basis $\{v\}$ of $V$ and identify $V$ with $A$
by sending $av$ to $a$.
Write
$$
\Lambda_1=T_pE_1\subset A,
\qquad \Lambda_2=T_pE_2\subset A.
$$
After multiplying the coordinate by a power of $p$, we may assume
$\Lambda_1\subseteq R$. We then have $G_K$-equivariant identifications
$$
E_i[p^\infty]\simeq A/\Lambda_i\qquad(i=1,2).
$$
Since all CM endomorphisms are defined over $K$, the action of
$G_K$ on $V$ is $A$-linear. Thus each $\sigma\in G_K$ sends $v$
to $u_\sigma v$ for some $u_\sigma\in A^\times$ and acts on
$V\simeq A$ by multiplication by $u_\sigma$.

By Lemma~\ref{adelic}(ii), for every sufficiently large $N$
and every $a\in R$, there exists $\sigma_a\in G_K$ that acts on $V$
by multiplication by $1+p^Na$ and acts trivially on $T_qE_i$
for all $q\ne p$ and $i=1,2$.

For $n>N$, let
$$
P_n=p^{-n}\bmod\Lambda_1\in E_1[p^\infty].
$$
Since $\Lambda_1\subseteq R$, the order of $P_n$ is at least $p^n$.
Only finitely many of these points can be defined over a fixed number
field. Thus (i) implies that, for every sufficiently large $n$, there
is $Q\in E_2(\overline K)_{\tors}$ such that $K(P_n)=K(Q)$.
Fix such an $n$ and a corresponding point $Q$. Write
$Q=Q_p+Q_{p'}$, where $Q_p$ has $p$-power order and $Q_{p'}$
has order prime to $p$. Under the identification
$E_2[p^\infty]\simeq A/\Lambda_2$, choose $\beta\in A$ such that
$$
Q_p=\beta\bmod\Lambda_2.
$$
The stabilizers of $P_n$ and $Q$ in $G_K$ agree, and each $\sigma_a$
fixes $Q_{p'}$. Hence, for every $a\in R$,
$$
a\in p^{n-N}\Lambda_1
\quad\Longleftrightarrow\quad
p^Na\beta\in\Lambda_2.
$$
Set
$$
J=\{a\in A:p^Na\beta\in\Lambda_2\}.
$$
Then
\begin{equation}\label{refine:intersection}
J\cap R=p^{n-N}\Lambda_1\subseteq pR.
\end{equation}
Lemma~\ref{refine:smallsubmodule} gives $J\subseteq pR$, so
$J=J\cap R$.

We claim that $\beta\in A^\times$. Otherwise multiplication by
$\beta$ on $A$ would have a nonzero kernel. After rescaling by a
power of $p$, this kernel contains an element of $R\setminus pR$.
That element belongs to $J$, contrary to
\eqref{refine:intersection}. Thus $\beta$ is invertible, and
$$
p^{-N}\beta^{-1}\Lambda_2
=J=p^{n-N}\Lambda_1.
$$
Consequently,
\begin{equation}\label{refine:homothetic}
\Lambda_2=p^n\beta\Lambda_1.
\end{equation}
Multiplication by $p^n\beta$ commutes with the $G_K$-action and
therefore gives a $G_K$-equivariant isomorphism
$$
u:T_pE_1\xrightarrow{\sim}T_pE_2.
$$
By the integral Hom form of Faltings' theorem~\cite{faltings2},
$$
\Hom_K(E_1,E_2)\otimes_{\Z}\Z_p
\simeq
\Hom_{G_K}(T_pE_1,T_pE_2).
$$
Put $H=\Hom_K(E_1,E_2)$. 
Let $f_1,\ldots,f_r$ be a $\Z$-basis of $H$.
By the above isomorphism, we may write
$$
u=\sum_{i=1}^r a_i T_p(f_i)
\qquad\text{with } a_i\in\Z_p.
$$
For each $i$, choose $b_i\in\Z$ such that $b_i\equiv a_i\pmod p$,
and set $\psi=\sum_{i=1}^r b_i f_i\in H$. Then
$$
T_p(\psi)=\sum_{i=1}^r b_i T_p(f_i)
\equiv \sum_{i=1}^r a_i T_p(f_i)=u\pmod p.
$$
Under the canonical identifications $T_pE_i/pT_pE_i\simeq E_i[p]$,
this congruence shows that $\psi:E_1[p]\to E_2[p]$ is an
isomorphism. Thus $\psi$ is a nonzero isogeny, and its kernel
contains no nonzero point of order $p$. Hence $p\nmid\deg\psi$.
Since this construction works for every prime $p$, no prime divides
the degrees of all $K$-isogenies from $E_1$ to $E_2$.
Therefore $c_K(E_1,E_2)=1$.
\end{proof}

\begin{rem}\label{refine:noconverse}
The hypothesis on the endomorphisms in
Theorem~\ref{refine:criterion} is automatic for non-CM curves.
The theorem does not assert that equality of the full collections
$\L(E_i;K)$ implies equality of the torsion-generated collections.
In the non-CM case, such an implication would force
$E_1\simeq_K E_2$ by Theorem~\ref{mainthm,nonCM}, giving a positive
answer to Question~\ref{question:nonCM} with an isomorphism
defined over $K$.
\end{rem}

\subsection{CM counterexamples}

To give a negative answer to Question~\ref{question1.9}, we use
nonprincipal ideals in the rings of integers of imaginary quadratic fields.
After constructing the pairs in
Theorem~\ref{mainthm,CM}, we use quadratic twists to obtain infinitely
many pairs over a fixed number field, as in
Corollary~\ref{mainthm,quadtwist}. We use standard CM theory;
see \cite{cox,silverman2}.

Let $q\geq5$ be a prime satisfying $q\equiv1\pmod4$.
Set $F_q=\Q(\sqrt{-q})$ and $\O_q=\Z[\sqrt{-q}]$.
Since $-q\equiv3\pmod4$, $\O_q$ is the ring of integers of $F_q$,
and $\mathrm{disc}(\O_q)=-4q$.
Consider the ideal
$$
\a_q=(2,1+\sqrt{-q})\subset\O_q.
$$

\begin{lem}\label{idealnorm2}
The ideal $\a_q$ has norm $2$ and is nonprincipal.
\end{lem}

\begin{proof}
We have $\O_q/\a_q\cong\F_2$, and therefore $\a_q$
has ideal norm $2$.

Suppose that $\a_q$ is principal.
Then there exists an element $u+v\sqrt{-q}\in\O_q$
of norm $2$. Thus
$$
u^2+qv^2=2.
$$
For $q\geq5$, this equation has no integer solutions.
Hence $\a_q$ is nonprincipal.
\end{proof}

Consider the complex elliptic curves
\begin{equation}\label{elleq}
E_{1,q}^{\an}=\C/\O_q,
\qquad
E_{2,q}^{\an}=\C/\a_q.
\end{equation}
Since $\a_q$ is an invertible ideal, both curves have
endomorphism ring $\O_q$.
If they were isomorphic, their corresponding lattices would
be homothetic, so $\a_q=u\O_q$ for some $u\in F_q^\times$.
This contradicts the fact that $\a_q$ is nonprincipal.
Hence $E_{1,q}^{\an}$ and $E_{2,q}^{\an}$ are not isomorphic.

Let $0\ne\beta\in\a_q$.
Multiplication by $\beta$ induces a homomorphism
$$
\phi_\beta:E_{1,q}^{\an}\longrightarrow E_{2,q}^{\an},
\qquad
z+\O_q\longmapsto\beta z+\a_q.
$$
It is an isogeny, and its degree is
$$
\deg\phi_\beta
=[\a_q:\beta\O_q]
=\frac{\mathrm{N}_{F_q/\Q}(\beta)}{\mathrm{N}(\a_q)}
=\frac{\mathrm{N}_{F_q/\Q}(\beta)}{2}.
$$
In particular, the elements $2$ and $1+\sqrt{-q}$ of $\a_q$
induce isogenies of degrees
$$
\deg\phi_2=2,
\qquad
\deg\phi_{1+\sqrt{-q}}=\frac{q+1}{2}.
$$
Since $q\equiv1\pmod4$, the integer $(q+1)/2$ is odd.
Hence these two degrees are relatively prime.

The $j$-invariants of the complex elliptic curves
$E_{1,q}^{\an}$ and $E_{2,q}^{\an}$ are algebraic.
Therefore these curves admit models over $\overline{\Q}$.
The homomorphisms induced by $2$ and $1+\sqrt{-q}$ are
algebraic and are defined over finite extensions of $\Q$.
Hence there exists a number field $K_q$ over which both
curves admit models and over which the two isogenies
$$
\phi_2:E_{1,q}\longrightarrow E_{2,q},
\qquad
\phi_{1+\sqrt{-q}}:E_{1,q}\longrightarrow E_{2,q}
$$
are defined.

We now prove Theorem~\ref{mainthm,CM}.
For the reader's convenience, we restate it here.

\begin{thm}[=Theorem~\ref{mainthm,CM}]\label{3.5}
For every prime $q\equiv1\pmod4$, there are a number field $K_q$
and elliptic curves $E_{1,q},E_{2,q}$ over $K_q$ such that
$$
\L(E_{1,q};K_q)=\L(E_{2,q};K_q),
\qquad
E_{1,q}\not\simeq_{\overline{K_q}}E_{2,q}.
$$
Moreover, as $q$ varies, all the curves $E_{1,q}$ and $E_{2,q}$
are pairwise non-isomorphic over $\overline{\Q}$.
\end{thm}

\begin{proof}
The two $K_q$-isogenies $\phi_2$ and $\phi_{1+\sqrt{-q}}$
have relatively prime degrees.
Hence Corollary~\ref{coprimedegreecor} gives
$$
\L(E_{1,q};K_q)=\L(E_{2,q};K_q).
$$
However, since their complex uniformizing lattices $\O_q$
and $\a_q$ are not homothetic, the curves are not
geometrically isomorphic.
If $q'$ is another prime satisfying $q'\equiv1\pmod4$,
then $\O_{q'}$ and $\O_q$ have different discriminants,
so their geometric endomorphism rings are not isomorphic.
Thus all the curves obtained as $q$ varies are pairwise
non-isomorphic over $\overline{\Q}$.
\end{proof}

\begin{exam}[The case of discriminant $-20$]\label{example}
For $q=5$, we have
$$
\O=\Z[\sqrt{-5}],
\qquad
\a=(2,1+\sqrt{-5}).
$$
The two isogenies $\phi_2$ and $\phi_{1+\sqrt{-5}}$ have degrees $2$ and $3$, while the two
curves are not geometrically isomorphic.
The two $j$-invariants are the roots of the Hilbert class polynomial
$$
H_{-20}(X)=X^2-1264000X-681472000,
$$
namely,
\begin{equation}\label{jeq}
j_{\pm}=632000\pm282880\sqrt{5}=(50\pm26\sqrt{5})^3.
\end{equation}
In particular, $j_+\ne j_-$.
\end{exam}

\begin{cor}\label{mainthm,quadtwist}
For some number field $K$, there are infinitely many pairs
$(E_{1}^{d},E_{2}^{d})$ of elliptic curves, pairwise distinct
up to $K$-isomorphism, such that
$$
\L(E_{1}^{d};K)=\L(E_{2}^{d};K),
\qquad
E_{1,d}\not\simeq_{\overline K}E_{2,d}.
$$
\end{cor}

\begin{proof}
Choose a number field $K$ over which the curves in
Example~\ref{example} and their isogenies of degrees $2$ and $3$
are defined. For each $d\in K^\times/(K^\times)^2$, twist both
curves by $d$. The isogenies twist to $K$-isogenies of the same
degrees, so Corollary~\ref{coprimedegreecor} gives equality of the
collections $\L$. The two $j$-invariants in \eqref{jeq} are distinct
and neither is $0$ nor $1728$. Thus the curves in each pair remain
geometrically non-isomorphic, and distinct square classes give
distinct $K$-isomorphism classes of twists. Since
$K^\times/(K^\times)^2$ is infinite, the result follows.
\end{proof}

\section{An application to rank growth}\label{ranksection}

Equality of the full collections also constrains rank growth, without
any assumption on complex multiplication.

\begin{prop}\label{sec2cor2}
Let $E_1,E_2/K$ be elliptic curves such that
$\L(E_1;K)=\L(E_2;K)$. For every finite extension $F/K$ and every
prime $p$, we have
$$
\rank E_1(L)>\rank E_1(F)
\quad\Longleftrightarrow\quad
\rank E_2(L)>\rank E_2(F)
$$
for all but finitely many cyclic extensions $L/F$ of degree $p$.
In particular, the quadratic twists of the two curves have positive
rank over $F$ for the same nontrivial quadratic characters of $G_F$, up to
finitely many exceptions.
\end{prop}

\begin{proof}
Fix $F$ and $p$. By Lemma~\ref{torsionbasechange}, the collections
remain equal over $F$. An extension $L/F$ of degree $p$ belongs
to $\L((E_i)_F;F)$ exactly when $E_i(L)\ne E_i(F)$, since there
are no proper intermediate fields. Hence
$$
E_1(L)\ne E_1(F)
\quad\Longleftrightarrow\quad
E_2(L)\ne E_2(F).
$$
It therefore suffices to show that, for a fixed elliptic curve
$E/F$, only finitely many cyclic degree-$p$ extensions give new
points without increasing the rank.

Suppose $E(L)\ne E(F)$ and $\rank E(L)=\rank E(F)$, and choose
$P\in E(L)\setminus E(F)$. Put
$$
R=\sum_{\sigma\in\Gal(L/F)}\sigma(P)\in E(F),
\qquad D=pP-R.
$$
The equality of ranks implies that a positive multiple of $P$
belongs to $E(F)$. If $nP\in E(F)$, then $nR=pnP$ and
$nD=0$, so $D$ is torsion. Since $D\in E(L)$, we have $[F(D):F] \le p$.
There are only finitely many such
torsion points by Northcott's theorem and the vanishing of their
canonical heights; see
\cite[Theorems~VIII.5.11 and~VIII.9.3(d,e)]{silverman1}.

Choose representatives $\mathcal R$ for the finite group
$E(F)/pE(F)$ and write $R=R_0+pQ$, with $R_0\in\mathcal R$ and
$Q\in E(F)$. Then $p(P-Q)=R_0+D$. The right-hand side belongs
to a fixed finite set, and multiplication by $p$ has finite
fibres. Thus $P-Q$ belongs to a fixed finite set. Since
$L=F(P)=F(P-Q)$, only finitely many such extensions occur.
Applying this to both curves proves the rank-growth assertion.

For a quadratic extension $L/F$ with associated character $\chi$,
the standard eigenspace decomposition gives
$\rank E_i(L)=\rank E_i(F)+\rank E_i^\chi(F)$, where
$E_i^\chi$ is the quadratic twist of $(E_i)_F$ by $\chi$.
This proves the final assertion.
\end{proof}

\section*{Acknowledgments}
The first author was supported by the National Research Foundation
of Korea (NRF) grant funded by the Korea government (MSIT)
(RS-2024-00342122). 
The second author was supported by the National Research Foundation
of Korea (NRF) grant funded by the Korea government (MSIT)
(RS-2025-23525445).

\section*{AI Declaration}
After making the observations in Theorem~\ref{degreecontentprop}
and Corollary~\ref{coprimedegreecor} and developing the idea of
using complex multiplication to construct counterexamples,
we used OpenAI's ChatGPT 5.6 Pro to find the counterexamples
in Theorem~\ref{3.5}. All other use of AI was limited to
grammar correction and language editing.


\end{document}